\documentstyle[12pt,epsf]{article}
\input epsf
\def\bee{\begin{equation}}
\def\ee{\end{equation}}

\begin{document}

{\large
  \centerline{\bf Repopulation: is it inevitably?}
}

\bigskip
\bigskip

\begin{center}
{\large K. Ma{\l}achowski$^1$,  M. Wolf$^2$\\}
\bigskip
\bigskip
{\small
\noindent $^1$ Lower Silesian Oncology Center\\
Pl. Hirszfelda 12\\
PL-53-413 Wroc{\l}aw, Poland, e-mail: malachowski.k@dco.com.pl\\

\bigskip
\bigskip
$^2$Institute of Theoretical Physics, University of Wroc{\l}aw\\
Pl.Maxa Borna 9, PL-50-204 Wroc{\l}aw, Poland, e-mail:
mwolf@ift.uni.wroc.pl\\
}
\end{center}

\bigskip

\begin{abstract}
 A mathematical model  of radiotherapy is proposed. The study used the
classical 24 hours way of fractionation with a weekend pause. We introduce the matrices of
``radiotherapy'' and ``growth''. We developed an equation of the fraction cell evolution, which we
solved numerically. The results indicate that the accelerated growth of cells occurs
due to the decrease of the fraction of slowly growing  cells and increase of the cells that are fast
growing.
\end{abstract}

\bigskip
\bigskip
\bigskip

Repopulation is a serious problem in cancer radiotherapy. The
growth of tumor in final stage of radiation can prevent the
healing of a patient. Many studies have shown the importance of
timing in radiotherapy \cite{Maciej}. If there is a delay in the treatment
caused e.g. by interrupts in radiation then this the increases time given to cells for accelerated
growth.  It was shown that breaks in radiation especially after fourth weeks of treatment
lead
to the worse results, approximately 2--4.8 percent growth per day of delay.

Why such a strange phenomenon occurs?  Tumor fights for its live and it  behaves according to
the Lenz rule: the lower number of cells the faster the growth.
According to Trott \cite{Trott} accelerated repopulation occurs when the number of cells in
tumor decreases
below 1000 cells. It can be explained in a few different ways.

Referring to a recent study there \cite{Gasinska}  are  three theories for explaining reasons for
repopulations:

1.  Fowler model, in which the author claims  that cancer volume
doubling time $T_d$ approaches potential time $T_{pot}$ as a
result of loosing cells.

2. Jones Model which proposes the following explanation:  the tumor
possesses subpopulations of cells growing with different
velocities (speeds); cells are dying equally but those dividing
faster gain advantage during breaks in radiation.

3.   Trott--Kummermehr model which can be called Dragon Theory. Like
a mediewal knight cutting dragon's head have met next two new
heads, here stem  cells switch from asymmetrical division to
symmetrical one. At each division from one stem cell two are
arising.

\section{Assumptions}

We assume as a basis the Jones model.Cancer tumor is heterogenic; it means that there exists
fractions of cells that  differ in terms of access to oxygen or nutrition or number of mutations.We assume there are three fractions of cells in the tumor, which
we will denote:

$ x_0$ --  small number of mutations and low growth velocity

$x_1$ -- intermediate number of  mutations and medium growth velocity

$x_2$ -- large number of mutations and fastest growth velocity

where these variables are normalized by

$$ \sum_{i=0}^2 x_i =1 $$

In individual fractions  there is a well determined number of
cell, where

$ y_i$ number of cells belonging to the $i$-th fraction

$$ \sum_{i=1}^3 y_i = N $$

and where $N$ is a total number of cells.

One of mutation factors is the radiation itself \cite{Gordon}. There is no reason
to  prevent such a phenomenon during radiation. Subsequently to the radiation
of tumor  after each consecutive dose,  the number of cells in
tumor will decrease and cells in each individual fraction will
undergo   mutations. Below is a new model of decreasing of number of cells in each fraction.

Fraction $x_i^0$ can be expressed by

\bee
x_i^0= y_i^0/N
\ee

Growth of tumor is a result of growth of individual fractions of tumor cells.
Each fraction $x_0^0,~x_1^0, ~x_2^0$ grows with its own velocity
$v_0,~v_1, ~v_2$. Time  of duplication of tumor $T_{d,i}$ determines the velocity
according to:
\bee
v_i=\frac{\ln(2)}{T_{d,i}}
\ee
where $T_{d,i}$ is volume doubling time for individual fraction tumor.There is  besides velocity the influence on the
tumor growth,on  the number  of cells in the particular fractions of the tumor.
Average velocity of the  tumor growth  can be described by means of the formula:
\bee
\Phi = \sum_{i=0}^2 v_i x_i
\ee

\section{Matrix of radiation}

The decrease of the tumor volume, i.e. waste of stem cells is
described by the linear -- quadratic formula:

\bee
N=N_0 e^{-\alpha d -\beta d^2}
\ee

where $N$ is the  number of survive cells that radiotherapy , and $N_0$ is the initial number of cells.
Another form of this formula is following:

\bee
 S=  e^{-\alpha d -\beta d^2}
\ee

where $S$ is fraction of surviving cells, $\alpha, ~\beta$ are coefficients,
$d$ is a radiation dose. Accordingly
\bee
S=N/N_0
\ee
or
\bee
S=\sum_{i} x_i^1~~~~~~~~~~~~~~1=\sum_i x_0^i
\ee
Here $x_i^0$ are initial fractions before radiotherapy and $x_i^1$ after first dose of
radiation.

\bee
\sum_i x_i^1= \sum_i x_i^0 e^{-\alpha d -\beta d^2}
\ee

When move from fractions to the number of cells then the equations take the form:
\begin{eqnarray}
y_0^1 & = & y_0^0  (e^{- (\alpha d + \beta d^2)})\\
y_1^1 & = & y_1^0 e^{- (\alpha d + \beta d^2)} \\
y_2^1 & = & y_2^0 e^{- (\alpha d + \beta d^2)}
\end{eqnarray}

In matrix notation we have:

\begin{equation}
\pmatrix{y_0^1 \cr y_1^1 \cr y_2^1 } = \pmatrix{e^{-\alpha d - \beta d^2} & 0 & 0 \cr
                                                0 & e^{-\alpha d - \beta d^2} & 0 \cr
                                                0 & 0 &  e^{-\alpha d - \beta d^2}}
\pmatrix{y_0^0 \cr y_1^0 \cr  y_2^0 }
\end{equation}

The matrix is diagonal and it shows that each fraction decreases according to the
linear-quadratic formula and there is no exchange of cells between individual
fractions.

This description suggests that all cells behave the same way and
are equally sensitive to the absorbed dose gained by the tumor.
Investigations show that cancer tumor does not possess uniform cells
the individual cells differ in access to the oxygen or nutritious means.
One of mutation factors is the ion radiation.  During radiation surviving  cells inherit improved
conditions of oxygenations and nutrition  and are undergoing rapid mutation. All these changes
lead both  to the decrease in the number of cells in individual fractions and also to the change
of the proportions of individual fractions.   We introduce coefficients $Q$ and $P$
to describe the probability that the cells from fractions $x_0^0, x_1^0$ will shift
to fraction $x_1^0$ and $x_2^0$, respectively.

Equation describing this process have the following form:
\begin{eqnarray}
y_0^{1} & = & y_0^0({e^{-\alpha d - \beta d^2}}  - Q) \\
y_1^{1} & = & y_1^{0}({e^{-\alpha d - \beta d^2}} - P) + Q y_0^0\\
y_2^{1} & = & y_2^{0}{e^{-\alpha d - \beta d^2}}+P y_1^0 \\
\end{eqnarray}

Situation after $n$ steps is described by the following equations:

\begin{eqnarray}
y_0^{(n)} & = &\left( {e^{-(\alpha d +\beta d^2)}} - Q)\right) ~~ y_0^{(n-1)} \\
y_1^{(n)} & = & Qy_0^{(n-1)} + \left({e^{-(\alpha d +\beta d^2)}}-P)\right) y_1^{(n-1)} \\
y_2^{(n)} & = & P(y_1^{(n-1)} + y_2^{(n-1)}{e^{-(\alpha d +\beta d^2)}}  \\
\end{eqnarray}

or in matrix notation:

\begin{equation}
\pmatrix{y_0^{(n)} \cr y_1^{(n)} \cr y_2^{(n)}} = \pmatrix{e^{-\alpha d - \beta d^2}-Q & 0
& 0 \cr
                                                Q & e^{-\alpha d - \beta d^2} -P & 0 \cr
                                                0 & P &  e^{-\alpha d - \beta
d^2}}^{\Large n}
\pmatrix{y_0^0 \cr  y_1^0 \cr y_2^0 }
\end{equation}

We will call the matrix  appearing above a radiation matrix and we will denote
 it $\mathcal{R}$:

\begin{equation}
\mathcal{R} = \pmatrix{e^{-\alpha d - \beta d^2}-Q & 0 & 0 \cr
                                                Q & e^{-\alpha d - \beta d^2} -P & 0 \cr
                                                0 & P &  e^{-\alpha d - \beta
d^2}}^{\Large n}
\end{equation}

It can be shown by induction that for each $n$ after summing up rows we obtain
following
equations:

\bee
N=N_0 ~ e^{-n(\alpha d +\beta d^2)}
\ee
It means that the matrix $\mathcal{R}$ describes the diminishment of tumor cells  after
 radiation
according to the linear  quadratic form. Additionally  it shows how individual fractions
change in time during radiation. Coefficients $Q$ and $P$ allows exchange of cells between
individual fractions.

Inserting here values for $\alpha,~ \beta, ~P,~Q$ the number of cells $N$ and
values of fractions $x_0,~x_1$ we can calculate the rate of decrease of the number of
cells in each fraction. From computer simulations it follows that the vector
$(x_1,x_2,x_3)$
tends to the equilibrium state $(0,0,x_2)$. From this we conclude, that radiation of the
tumor leads to the selection of cells which are the most mutated and which grow with the
 largest speed.
In the limit of  large $n$ this equation has
the form:
$$ \Phi = v_2 x_2^{(n)}$$

It means that speed of tumor growth is larger when the tumor diminishes, and this effect cannot
be avoided.

\section{The growth matrix}

We describe the rate of tumor growth according to the Sole  \cite{hiszpanie}.
We make use of the results of the paper \cite{hiszpanie} in which
it was shown that the equations for growth of the cell fractions can be written as:

\begin{eqnarray*}
\frac{dx_0}{dt} & = & v_0x_0(1-Q') -x_0\Phi(x_0,x_1,x_2)\\
\frac{dx_1}{dt}& = & v_1x_1(1-P') + v_0x_0Q'  -x_1\Phi(x_0,x_1,x_2)\\
\frac{dx_2}{dt} & = & v_2x_2+ v_1 x_1 P'  -x_2\Phi(x_0,x_1,x_2)
\end{eqnarray*}

or in the matrix notation:

\bee
\dot{\overrightarrow{x}} = {\mathcal{M}}\overrightarrow{x}
\ee
where $\dot{\overrightarrow{x}}$ is the time derivative of the $\overrightarrow{x} =
(x_1, x_2, x_3)^T$, ${\mathbf{1}}$ is the identity matrix and mixing matrix $\mathcal{M}$
is given by
\bee
{\mathcal{M}} = \pmatrix{f_0(1-Q') - \Phi(x_0,x_1,x_2) & 0 & 0 \cr
                                               f_0 Q' & f_1(1-P') - \Phi(x_0,x_1,x_2) & 0
\cr
                                                0 & f_1 P'  &  f_2 - \Phi(x_0,x_1,x_2)}.
\ee

and we denote this matrix $\mathcal{M}$ as the growth matrix. Equations of Sole  describe
how ratios of fractions are changing during radiotherapy. However growth of the cell number
in each fraction which occurs during the pauses in radiotherapy we obtain through the following
procedure: The result of the Sole equation expressed in fractions $x_i$ we convert to the
integer valued number of cells $y_i$ (see below). Obtained number of cells we multiply by
factor $e^{ln(2)*v_i} = 2^{v_i}$ and next we pass from the number of cells back to fractions:
\bee
\pmatrix{y'_0 \cr y'_1 \cr y'_2 } = \pmatrix{2^{v_0} & 0 & 0 \cr
                                                0 & 2^{v_1} & 0 \cr
                                                0 & 0 &  2^{v_2}}
\pmatrix{y_0  \cr  y_1 \cr y_2}
\ee
Here $y'_i$  is the number of cells after division during pause between consecutive
pulses of radiotherapy. Let
\bee
{\mathcal{D}} =\pmatrix{2^{v_0} & 0 & 0 \cr
                                                0 & 2^{v_1} & 0 \cr
                                                0 & 0 &  2^{v_2}}
\ee
denote the division matrix.
The growth matrix we construct from mixing and division in the following way:
\bee
{\mathcal{G}}(y_i, x_i) = {\mathcal{D}}(y_i) {\mathcal{M}}(x_i)
\ee

We introduce the dependence of $v_2$  on $v_1$:
\bee
v_2=a v_1 \psi(\theta,Q,P,Q',P',d,n).
\ee
where
\bee
\psi(\theta,Q,P,Q',P',d,n)=e^{\Theta - n(\sqrt{Q^2+P^2}d+\sqrt{Q'^2+P'^2}d^2}
\ee
This form is for radiation period, during weekend the form is different:
\bee
\psi(\theta,Q',P',d,n)=e^{\Theta - n(\sqrt{Q'^2+P'^2}}
\ee
It contains a threshold after crossing this threshold velocity $v_2$ decreases, in
 accordance velocity
$\Phi$ also diminishes. In the paper by Sole at al \cite{hiszpanie} it is shown for
 which parameters
$Q^\prime,
P^\prime, M^\prime$ the fraction $x_2$ exists.

\section{The radiation}

The act of radiation consists of:

a. time of radiation --- pulse radiation

b. time  between consecutive exposures

These equations describe growth of tumor from the moment of ending of radiation till
the next exposures. In classical fraction this time is 24 hours. The radiation lasts for a very
short period of time: a few minutes, in comparison to 24 hour waiting period between
consecutive fractions. Exposition after Wheldon \cite{Wheldon}  we can call pulse
radiotherapy. In the remaining time, tumor cells repair damage from radiotherapy and they divide. The tumor growth
appears.

In classical radiotherapy we radiate once a day during the 4-7 weeks, depending on the
 radiation  dose quantity.
Duration of radiation is very short (a few minutes), remaining time is spent on the repair
of post-radiation damage. Symbolically we can demonstrate this in
the following form:
\begin{eqnarray*}
\nonumber \\
\label{iteracje}
(growth~ radiation)^5 growth^2 \ldots growth^2 (growth ~ radiation )^5
|N_0> =\\
 ((growth ~ radiation )^5) (growth^2 (growth ~ radiation )^5)^{n-1}|N_0>
\end{eqnarray*}
where $n$ is the number of weeks of radiotherapy  and $N_0$ is initial number of tumor cells,
i.e. $|N_0> = |initial~ number~of~tumor~cells~>$

Radiation is represented by the matrix:

\begin{equation}
\mathcal{R} = \pmatrix{e^{-\alpha d - \beta d^2}-Q & 0 & 0 \cr
                                                Q & e^{-\alpha d - \beta d^2} -P & 0 \cr
                                                0 & P &  e^{-\alpha d - \beta d^2}}
\end{equation}

and

\bee
{\mathcal{M}} = \pmatrix{f_0(1-Q') - \Phi(x_0,x_1,x_2) & 0 & 0 \cr
                                               f_0 Q' & f_1(1-P') - \Phi(x_0,x_1,x_2) & 0
\cr
                                                0 & f_1 P'  &  f_2 - \Phi(x_0,x_1,x_2)}.
\ee

The radiation process can be described by the following equation:
\bee y^{(n+1)} =  {\mathcal{R}} y^n \ee This equation describes
the diminishing of the number of cells in fractions after a pulse
of radiation. Next we can  write the equation which describe the
rise of tumor until next exposition. First we describe the change
in proportion in fractions in tumor cells. To this aim we change
the variable from $y$ to $x$. We assume that the proportion of
$y(t+{\rm next~ day})/y(t)$ is the same as normalized variables:

\bee \frac{y(t+{\rm next~ day})}{y(t)} = \frac{x(t+{\rm next
~day})}{x(t)} \ee
The Sole equations describe change of proportion
during growth of tumor between expositions: \bee
\frac{dx^{n+1}}{dt} = {\mathcal{M}}x^n \ee Differential equations
were solved on the interval of one  working day and on the
interval of three days during  weekend.We change solutions of this
equation  again into the number of cells.Then we calculate how
many new cells arise during the pause between radiations. We again
pass from $x$ back to $y$ as before and we use following equation:
\bee \overrightarrow{y}^{n+1} =  {\mathcal{V}}
\overrightarrow{y}^{n+1} \ee We repeat this procedure during
prescribed cure time expressed in weeks.

\section{Results}

We performed numerical calculations for two different sets of
parameters: in the first case the probability coefficients of
$Q,P,Q^\prime,P^\prime$ were zero, in the second they were
different from zero  and were
$Q=0.0005,P=0.0005,Q^\prime-0.1,P^\prime=0.1$. We took the
parameters $\alpha=0.2, \beta=0.02, d=2Gy, n=30, v_0=0.01,
v_1=0.016$. $v_2$ was calculated from the formula \ref{wzor2} for
$m=5$ and threshold $\theta=0.005$. At $Q,P,Q^\prime,P^\prime$
equal zero velocity of growth rising of tumor was initially
slower, and next faster. The conclusion is simply: the fact of
existence of the population of cells of with different growth
rates is sufficient for the tumor to grow faster as the corollary
of the radiation and intervals between the fraction radiotherapy.
At this point the faster   growing cells are gaining the
population dominance  over the slower growing cells and finally
lead to the accelerated growth tumor. After introduction the
$P,Q,Q^\prime,P^\prime$ different from zero the velocity of the
tumor growth was also larger. According to this model
probabilities $Q,~P,~Q',~ P'$ are responsible for the velocity
change via the change of the cell population $x$ from $x_0$ and
$x_1$ to $x_2$. Fraction $x_2$ is the fastest growing.
Biologically it can be explained that fraction $x_2$ gains the
best condition for growth due to the improvement of oxygenation
and better nutrition. Additionally fraction $x_2$ is built from
the most undifferentiated cells and most mutated cells. The factor
that causes the decrease fractions $x_0$ and $x_1$ and increase of
the fraction $x_2$ is the radiotherapy which leads to the death of
cells and simultaneously by causing the shift of  mutation from
$x_0$ and $x_1$ to $x_2$. In this way additional radiation
influence leads to some differences between situation when the
ionization energy is the cause of destroying only cells in
fractions ($Q = 0,P = 0,Q^\prime = 0,P^\prime = 0$) and leads to
the shift of cells from one fraction to another ($Q \neq 0,P \neq
0,Q^\prime \neq  0,P^\prime \neq 0$). These differences can be
seen after subtraction of two graph, see Figs. 1 and 2. In this
way we can see how the change of the oxygenation, better nutrition
and mutation influence  growth velocity. According to the
assumption about the existence of the threshold $\Theta$ the
influence of the coefficients ($Q \neq 0,P \neq  0,Q^\prime \neq
0,P^\prime \neq 0$)is  as following: In the first stage we see
that they accelerate tumor growth and next they cause the slowing
down of the tumor growth.  On this example we see a new mechanism
leading to the death of the tumor. Namely the existence of the
$\Theta$ threshold, Namely the existence of  the theta threshold,
causes that growth mutacji not only shifts cells to faster growing
fraction but also crossing sufficient threshold causes slower
growth of tumor.


\section{Conclusions}

Our model explains accelerated growth tumor according to Jones
model. The existence of cells fraction of different growth
velocity is the cause of the accelerated growth tumor. During the
radiation  therapy because of the interrupts, cells have time to
take advantage of the differences in the velocities of growth and
increase the number of cells in fastest growing fraction. In the
course of radiation the living conditions of cells are changing:
the oxygenation and nutrition is better and additionally mutations
appear. Together these factors causes the change of the number of
cells in the particular fractions. We have described these changes
by the coefficients $P$ and $Q$ responsible for the changes during
the radiation pulse and coefficients $Q^\prime ,P^\prime $
responsible for the changes occurring between different fractions.
These  coefficients are modifying the tumor growth: in the
beginning they accelerate and after crossing the threshold
$\Theta$ they slow down. The existence of the threshold $\Theta$
could explain the  benefits from the simultaneous
radiochemotherapy.  Chemotherapy has the mutagen function. This
mechanism in conjunction with radiotherapy facilitates the
crossing of the threshold $\Theta$ after which the tumor growth is
slowing down. From the model it follows that natural state of the
tumor is the state described by the vector $(0,0,1)$, it means
that it tends to the fastest fraction. The fastest and least
differentiated fraction gains the crucial dominance. It seems to
be in accordance with the clinical experience that often  the
revival of the tumor is more malicious and less differentiated.

Model that we introduced is the enlargement of the linear--quadratic model.
When we take into account only total number of cells it describes the diminishing
of the cells tumor exactly the same way as the linear--quadratic model. However it
allows to see how the numbers of cells in each fractions  change during the radiation
and which influence on the velocity of the tumor growth is the appearance of the
threshold $\Theta$ mutation.

It is possible to apply our model to arbitrary doses, time of radiation and breaks in
radiations time.\\

\vfill

\newpage

{\small
\vskip 0.4cm
\begin{center}
{\sf Table I}\\
\bigskip
\begin{tabular}{|c|c|c|c|c|} \hline
day  &  fraction 1  & fraction 2  & fraction 3  &  velocity  \\ \hline
1 & 371270035 &  210386353 & 37127004 &  0.000000000 \\ \hline
1 & 371476229 &  212652573 & 41821898 &  0.016515136 \\ \hline
2 & 229863321 &  131585880 & 25878696 &  0.016515136 \\ \hline
2 & 229878511 &  132938244 & 29136931 &  0.017018336 \\ \hline
3 & 142245005 &  82259977 & 18029449 &  0.017018336 \\ \hline
3 & 142177843 &  83060670 & 20288506 &  0.017571932 \\ \hline
4 & 87977288 &  51396563 & 12554190 &  0.017571932 \\ \hline
4 & 87883740 &  51866147 & 14118852 &  0.018180004 \\ \hline
5 & 54380999 &  32093910 & 8736511 &  0.018180004 \\ \hline
5 & 54287929 &  32366123 & 9818991 &  0.018846713 \\ \hline
6 & 53706425 &  32541035 & 11729058 &  0.019269721 \\ \hline
7 & 52629285 &  32602887 & 14884639 &  0.019715495 \\ \hline
8 & 32566127 &  20174125 & 9210367 &  0.019715495 \\ \hline
8 & 32380128 &  20263716 & 10310082 &  0.022962204 \\ \hline
9 & 20036285 &  12538851 & 6379707 &  0.022962204 \\ \hline
9 & 19900781 &  12581216 & 7133890 &  0.024044202 \\ \hline
10 & 12314273 &  7785047 & 4414332 &  0.024044202 \\ \hline
10 & 12217058 &  7802451 & 4930553 &  0.025208646 \\ \hline
11 & 7559712 &  4828027 & 3050944 &  0.025208646 \\ \hline
11 & 7490861 &  4832904 & 3403560 &  0.026457045 \\ \hline
12 & 4635221 &  2990521 & 2106067 &  0.026457045 \\ \hline
12 & 4587001 &  2989628 & 2346407 &  0.027789957 \\ \hline
13 & 4457412 &  2952492 & 2753154 &  0.028616727 \\ \hline
14 & 4249783 &  2878036 & 3399291 &  0.029472506 \\ \hline
15 & 2629695 &  1780881 & 2103425 &  0.029472506 \\ \hline
15 & 2583253 &  1767292 & 2326276 &  0.035259776 \\ \hline
16 & 1598474 &  1093571 & 1439461 &  0.035259776 \\ \hline
16 & 1567496 &  1083328 & 1589182 &  0.037028960 \\ \hline
17 & 969941 &  670345 & 983359 &  0.037028960 \\ \hline
17 & 949425 &  662867 & 1083679 &  0.038852571 \\ \hline
18 & 587489 &  410171 & 670562 &  0.038852571 \\ \hline
18 & 573996 &  404843 & 737601 &  0.040720840 \\ \hline
19 & 355179 &  250510 & 456415 &  0.040720840 \\ \hline
19 & 346366 &  246788 & 501095 &  0.042622914 \\ \hline
20 & 326744 &  236600 & 570777 &  0.043757469 \\ \hline
21 & 297833 &  220498 & 673761 &  0.044897230 \\ \hline
22 & 184294 &  136440 & 416912 &  0.044897230 \\ \hline
22 & 178080 &  133186 & 453545 &  0.051789508 \\ \hline
23 & 110193 &  82413 & 280646 &  0.051789508 \\ \hline
23 & 106278 &  80297 & 304735 &  0.053649548 \\ \hline
\end{tabular}
\end{center}
\vskip 0.4cm

}

\newpage

\newpage

{\small
\vskip 0.4cm
\begin{center}
\bigskip
\begin{tabular}{|c|c|c|c|c|} \hline
day  &  fraction 1  & fraction 2  & fraction 3  &  velocity  \\ \hline
24 & 65763 &  49686 & 188565 &  0.053649548 \\ \hline
24 & 63311 &  48322 & 204376 &  0.055462791 \\ \hline
25 & 39176 &  29901 & 126465 &  0.055462791 \\ \hline
25 & 37648 &  29028 & 136826 &  0.057219673 \\ \hline
26 & 23296 &  17962 & 84665 &  0.057219673 \\ \hline
26 & 22349 &  17408 & 91445 &  0.058911922 \\ \hline
27 & 20408 &  16155 & 100824 &  0.059878183 \\ \hline
28 & 17720 &  14341 & 113371 &  0.060818297 \\ \hline
29 & 10965 &  8874 & 70152 &  0.060818297 \\ \hline
29 & 10444 &  8539 & 75231 &  0.065925766 \\ \hline
30 & 6463 &  5284 & 46552 &  0.065925766 \\ \hline
30 & 6148 &  5078 & 49860 &  0.067155054 \\ \hline
31 & 3805 &  3142 & 30853 &  0.067155054 \\ \hline
31 & 3615 &  3017 & 33006 &  0.068300215 \\ \hline
32 & 2237 &  1867 & 20424 &  0.068300215 \\ \hline
32 & 2123 &  1790 & 21826 &  0.069362910 \\ \hline
33 & 1314 &  1108 & 13506 &  0.069362910 \\ \hline
33 & 1246 &  1061 & 14418 &  0.070345563 \\ \hline
34 & 1112 &  962 & 15538 &  0.070889546 \\ \hline
35 & 934 &  826 & 16893 &  0.071407364 \\ \hline
36 & 578 &  511 & 10453 &  0.071407364 \\ \hline
36 & 546 &  488 & 11117 &  0.074038283 \\ \hline
37 & 338 &  302 & 6879 &  0.074038283 \\ \hline
37 & 319 &  288 & 7312 &  0.074629103 \\ \hline
38 & 197 &  178 & 4524 &  0.074629103 \\ \hline
38 & 186 &  170 & 4806 &  0.075165861 \\ \hline
39 & 115 &  105 & 2974 &  0.075165861 \\ \hline
39 & 109 &  100 & 3158 &  0.075652632 \\ \hline
40 & 67 &  62 & 1954 &  0.075652632 \\ \hline
40 & 63 &  59 & 2074 &  0.076093354 \\ \hline
41 & 56 &  53 & 2209 &  0.076333559 \\ \hline
42 & 46 &  45 & 2362 &  0.076559766 \\ \hline
43 & 29 &  28 & 1461 &  0.076559766 \\ \hline
43 & 27 &  26 & 1548 &  0.077673641 \\ \hline
44 & 17 &  16 & 958 &  0.077673641 \\ \hline
44 & 16 &  15 & 1015 &  0.077915957 \\ \hline
45 & 10 &  10 & 628 &  0.077915957 \\ \hline
45 & 9 &  9 & 665 &  0.078133711 \\ \hline
46 & 6 &  6 & 412 &  0.078133711 \\ \hline
46 & 5 &  5 & 436 &  0.078329254 \\ \hline
47 & 3 &  3 & 270 &  0.078329254 \\ \hline
47 & 3 &  3 & 285 &  0.078504740 \\ \hline
48 & 3 &  3 & 303 &  0.078599769 \\ \hline
\end{tabular}
\end{center}

\newpage

\epsfverbosetrue

\begin{figure}[pht]
\vspace{-2cm}
\centering
\begin{minipage}{15.8cm}
\hspace{3.5cm}

\includegraphics{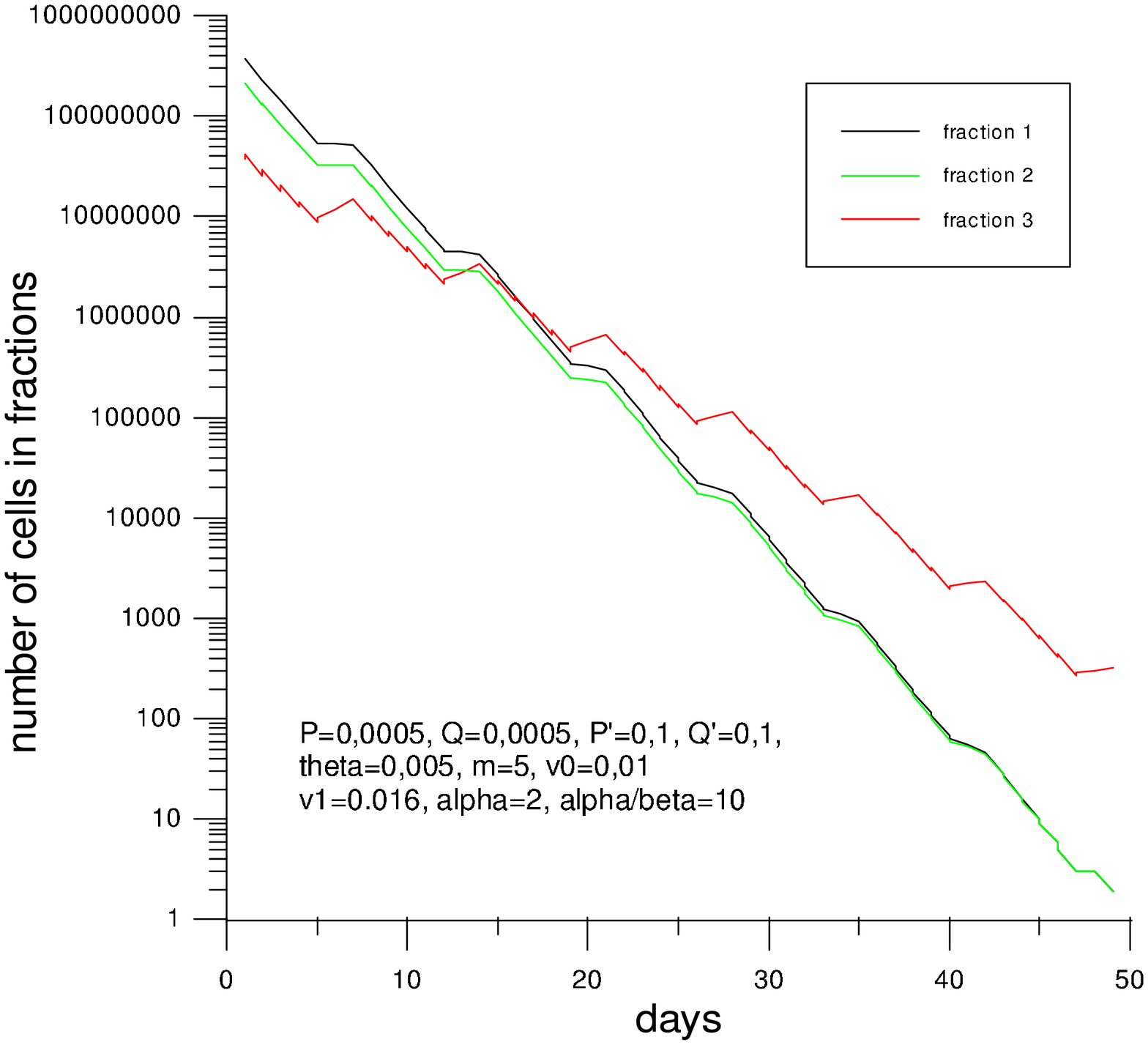}
\vspace{11.0cm}
\caption{Plot of velocity of the   tumor growth as a function of time (days)
 of radiotherapy. Here $Q, P, Q^\prime, P'$ are  zero.}
\label{comparison}
\end{minipage}

\begin{minipage}{15.8cm}
\hspace{0.0cm}
\includegraphics{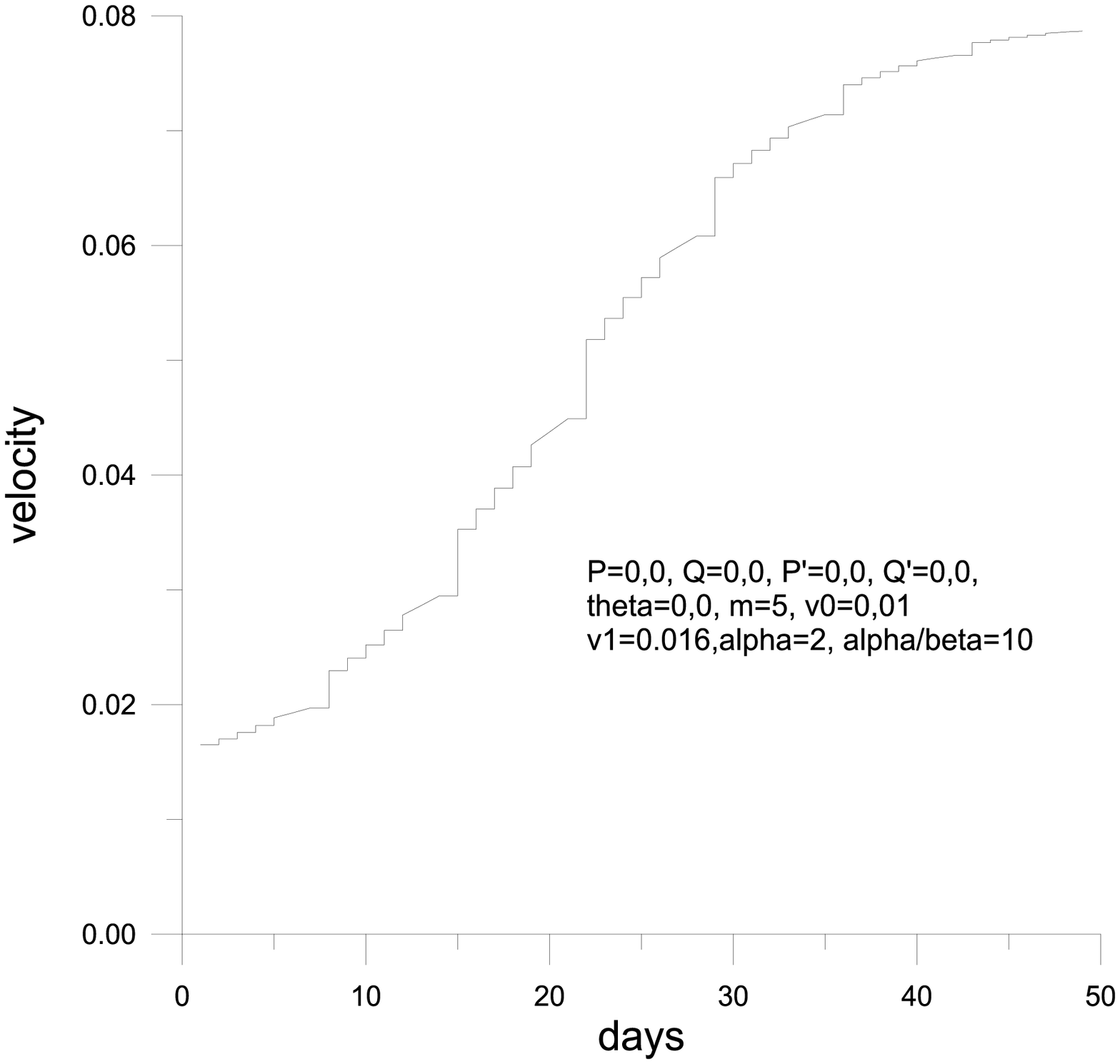}
\vspace{11cm}
\caption{Plot of velocity of the  tumor   growth as a function of time (days)
 of radiotherapy. Here $Q, P, Q^\prime, P'$ are different from zero.}.
\label{distr}
\end{minipage}
\end{figure}
\bigskip
\bigskip

\vfill

\newpage

\begin{figure}[pht]
\vspace{-4cm}
\centering
\hspace{2.5cm}

\includegraphics{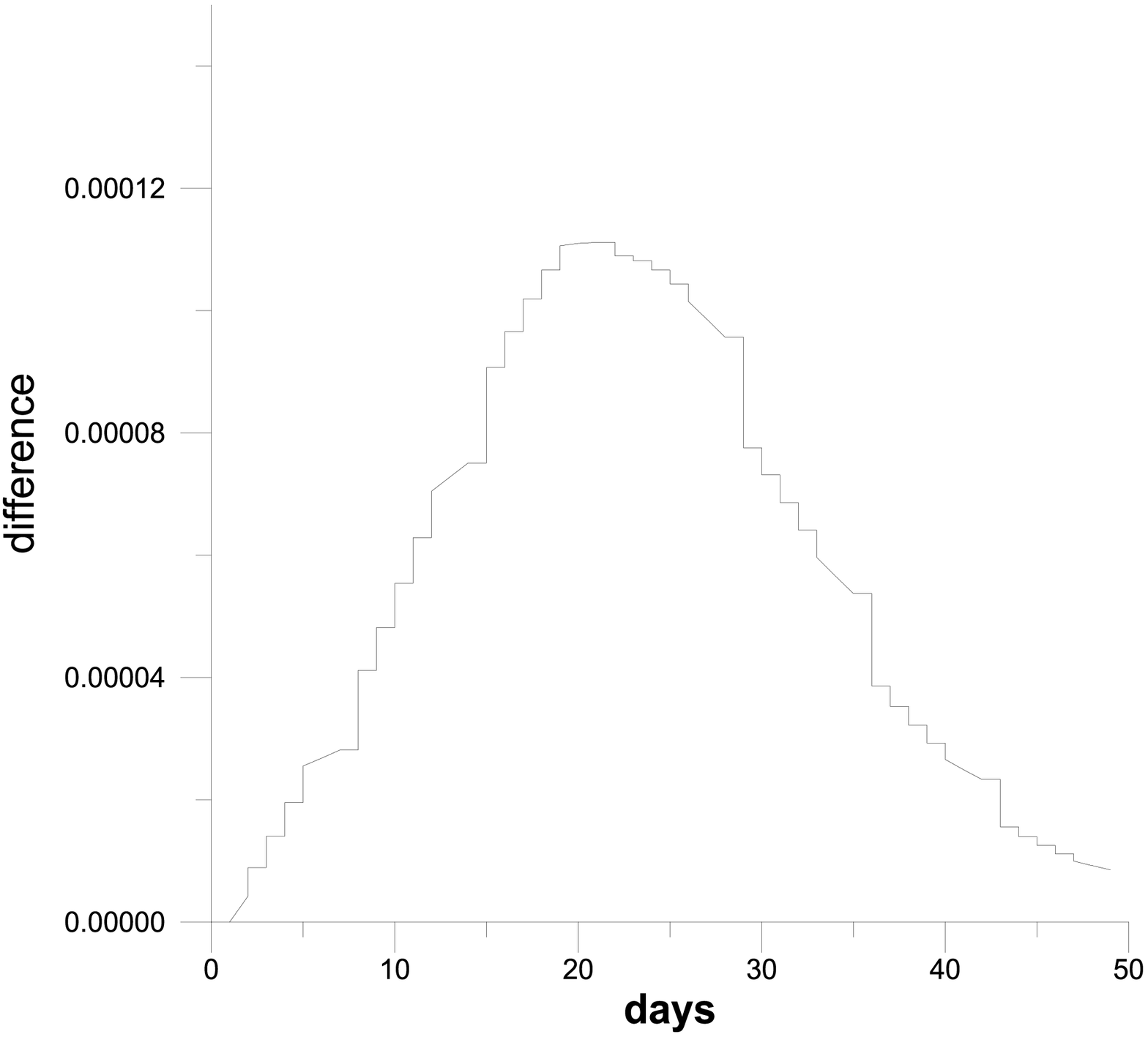}
\vspace{11.0cm}
\caption{Plot of velocity for zero valued parameters subtracted from velocity for parameters
different from zero.}
\label{comparison}
\end{figure}

\vfill

\newpage

\end{document}